\documentclass{amsart}

\usepackage{amssymb}
\usepackage{graphicx}

\newtheorem{thm}{Theorem}[section]

\newtheorem{pro}[thm]{Proposition}
\newtheorem{lem}[thm]{Lemma}

\newtheorem{cor}[thm]{Corollary}

\begin{document}

\title{A polynomial-time solution to the reducibility problem}
\author[K.~H.~Ko]{Ki Hyoung Ko}
\email{\{knot,leejw\}@knot.kaist.ac.kr}
\author[J.~W.~Lee]{Jang Won Lee}
%\email{leejw@knot.kaist.ac.kr}
\address{Department of Mathematics, Korea Advanced Institute
of Science and Technology, Daejeon, 305-701, Korea}
%%\thanks{}
%\author{Ki Hyoung Ko \and Jang Won Lee}
%\institute{Department of
%Mathematics, Korea Advanced Institute of Science and Technology,
%Daejeon, 305-701, Korea}%%\keywords{}
%\offprints{}
%\date{Received: date / Revised version: date}
\maketitle

\begin{abstract}
We propose an algorithm for deciding whether a given braid is
pseudo-Anosov, reducible, or periodic. The algorithm is based on
Garside's weighted decomposition and is polynomial-time in the
word-length of an input braid. Moreover, a reduction system of
circles can be found completely if the input is a certain type
of reducible braids.
\end{abstract}

\section{Preliminaries and introduction}
As a homeomorphism of a 2-dimensional disk that preserves $n$
distinct interior points and fixes the boundary of the disk, an
$n$-braid $x$ is isotopic to one of the following three dynamic
types known as the Nielson-Thurston classification \cite{T88}: (i)
{\em periodic} if $x^p$ is the identity for some nonnegative
integer $p$; (ii) {\em reducible} if $x$ preserves a set of
disjointly embedded circles; (iii) {\em pseudo-Anosov} if neither
(i) nor (ii). Obviously dynamic types are invariant under
conjugation and taking a power. A set of disjointly embedded
essential circles preserved by a reducible braid is called a {\em
reduction system}. Suppose the  distinct points lie on an axis. An
essential circle, i.e. separating $n$ distinct points, is {\em
standard} if it intersects the axis exactly twice. A reduction
system is {\em standard} if each circle in the system is standard.
Up to conjugacy, every reducible braid has a standard reduction
system. Standard reduction systems are especially nice in the
sense that they can be recognized in polynomial time.

Recently, some evidence that the conjugacy problem could be easy
for pseudo-Anosov braids has been found \cite{KL,BGG1}. In
addition, the conjugacy problem for periodic braids is trivial
once they are recognized. It is therefore important to know the
dynamic type of a given braid to solve its conjugacy problem. If a
braid is reducible, it is also important to know how the braid is
decomposed into pseudo-Anosov braids or periodic braids. Thus the
{\em reducibility problem} comes in two flavors depending on what
is asked. Given an arbitrary braid, we may ask to determine its
dynamic type or to find a reduction system if it is reducible. The
latter problem will be called the {\em reduction problem} to
distinguish them.

Our approach will be based on a weighted decomposition of braid
words, invented by Garside \cite{G}, improved by Thurston \cite{T}
and El-Rifai and Morton \cite{EM}. We briefly review the idea
together with necessary notations. The Artin presentation of the
group $B_n$ of $n$-strand braids has $n-1$ generators $\sigma_1,
\cdots, \sigma_{n-1}$ and two types of defining relations:
$\sigma_j\sigma_i =\sigma_i\sigma_j$ for $|i-j| > 1$ and
$\sigma_i\sigma_j\sigma_i = \sigma_j\sigma_i\sigma_j$ for $|i-j|=
1$. The monoid given by the same presentation is denoted by
$B_n^+$ whose elements will be called {\em positive braids}.

A partial order $\prec$ on $B_n^+$ can be given by saying $x\prec
y$ for $x, y\in B_n^+$ if $x$ is a ({\em left}) {\em subword} of
$y$, that is, $xz=y$ for some $z\in B_n^+$. Given $x,y\in B_n^+$,
the ({\em left}) {\em join} $x\vee y$ of $x$ and $y$ is the
minimal element with respect to $\prec$ among all $z$'s satisfying
that $x\prec z$ and $y\prec z$, and the ({\em left}) {\em meet}
$x\wedge y$ of $x$ and $y$ is the maximal element with respect to
$\prec$ among all $z$'s satisfying that $z\prec x$ and $z\prec y$.
Even though ``left" is our default choice, we sometimes need the
corresponding right versions: the partial order $\prec_R$ of being
a {\em right subword}, the {\em right join} $\vee_R$, and the {\em
right meet} $\wedge_R$. For example, $x\prec_R y$ if $zx=y$ for
some $z\in B_n^+$.

The {\em fundamental braid} $\Delta=(\sigma_1\cdots
\sigma_{n-1})(\sigma_1\cdots \sigma_{n-2})\cdots
(\sigma_1\sigma_2)\sigma_1$ plays an important role in the study
of $B_n$. Since it represents a half twist as a geometric braid,
$x\Delta =\Delta\tau(x)$ for any braid $x$ where $\tau$ denotes
the involution of $B_n$ sending $\sigma_i$ to $\sigma_{n-i}$. It
also has the property that $\sigma_i\prec\Delta$ for each
$i=1,\cdots,n-1$. Since the symmetric group $\Sigma_n$ is obtained
from $B_n$ by adding the relations $\sigma_i^2=1$, there is a
quotient homomorphism $q:B_n\to \Sigma_n$. For $S_n=\{x\in
B_n^+\mid x\prec\Delta\}$, the restriction $q:S_n\to\Sigma_n$
becomes a 1:1 correspondence and an element in $S_n$ is called a
{\em permutation braid}.

A product $ab$ of a permutation braid $a$ and a positive braid $b$
is ({\em left}) {\em weighted}, written $a\lceil b$, if $a^*
\wedge b=e$ where $e$ denotes the empty word and
$a^*=a^{-1}\Delta$ is the {\em right complement} of $a$. Each
braid $x\in B_n$ can be uniquely written as
$$x=\Delta^{u}x_1x_2\cdots x_{k}$$
where for each $i=1,\cdots,k$, $x_i\in S_n\setminus \{e,\Delta\}$
and $x_i\lceil x_{i+1}$. This decomposition is called the ({\em
left}) {\em weighted form} of $x$ \cite{G,T,EM}. Sometimes the
first and the last factors in a weighted form are called the {\em
head} and the {\em tail}, denoted by $H(x)$ and $T(x)$,
respectively. The weighted form provides a solution to the word
problem in $B_n$ and the integers $u$, $u+k$ and $k$ are
well-defined and are called the {\em infimum}, the {\em supremum}
and the {\em canonical length} of $x$, denoted by $\inf(x)$,
$\sup(x)$ and $\ell(x)$, respectively.

Given $x=\Delta^{u}x_1x_2\cdots x_{k}$ in its weighted form, there
are two useful conjugations of $x$ called the {\em cycling}
$\mathbf c(x)$ and the {\em decycling} $\mathbf d(x)$ defined as
follows: $${\mathbf c}(x)=\Delta^{u}x_2\cdots x_{k}\tau^{u}(x_1)=
\tau^{u}(H(x)^{-1})x\tau^{u}(H(x)),$$
$${\mathbf d}(x)=\Delta^{u}\tau^{u}(x_{k})x_1\cdots x_{k-1}=
T(x)xT(x)^{-1}.$$

%Since $\tau(\mathbf{d}(x))=\mathbf {c}(x^{-1})^{-1}$,
%$\mathbf{c}^{i}(x^{-1})=\tau^{i}(\mathbf{d}^{i}(x))^{-1}$ for
%$i\geq 1$.

A braid $x=\Delta^u x_1\cdots x_k$ in its weighted form is {\em
(left) $i$-rigid} for $1\le i\le k=\ell(x)$ if the first $i$
factors are identical in the weighted forms of $x_1\cdots x_k$ and
$x_1\cdots x_k\tau^u(x_1)$, that is, $x_1\cdots x_i=y_1\cdots y_i$
where $y_1\cdots y_{k}y_{k+1}$ is the weighted form of $x_1\cdots
x_k\tau^u(x_1)$ and $y_{k+1}$ may be empty. If a braid $x$ is
$\ell(x)$-rigid, we simply say $x$ is {\em (left) rigid} and this
is equivalent to the fact that $x_{k}\lceil\tau^{u}(x_1)$. We can
also consider the corresponding right version.

%A braid $x$ is {\em rigid} if its weighted form
%$x=\Delta^{u}x_1x_2\cdots x_{k}$ satisfies
%$x_{k}\lceil\tau^{u}(x_1)$. When $\ell(x)=1$,
%$x_1\lceil\tau^{u}(x_1)$ is required to be rigid. In
%some articles, the term ``rigid" has been used for this situation
%but we think our term is more suggestive.

Let $\inf_c(x)$ and $\sup_c(x)$ respectively denote the maximal
infimum and the minimal supremum of all braids in the conjugacy
class $C(x)$ of $x$. A typical solution to the conjugacy problem
in the braid group $B_n$ is to generate a finite set uniquely
determined by a conjugacy class. Historically the following four
finite subsets of the conjugacy class $C(x)$ of $x\in B_n$ have
been used in this purpose:

The {\em summit set}
$$SS(x)=\{y\in C(x)\mid \inf(y)=\mbox{$\inf_c$}(x)\}$$
was used by Garside in \cite{G} to solve the conjugacy problem in
$B_n$ for the first time. The {\em super summit set}
$$SSS(x)=\{y\in C(x)\mid \inf(y)=\mbox{$\inf_c$}(x)\mbox{ and }\sup(y)=\mbox{$\sup_c$}(x)\}$$
was used by El-Rifai and Morton in \cite{EM} to improve Garside's
solution. The {\em reduced super summit set}
$$RSSS(x)=\{y\in C(x)\mid \mathbf c^M(y)=y=\mathbf d^N(y)
\mbox{ for some positive integers }M,N\}$$ was used by Lee in his
Ph.D. thesis \cite{L00} to give a polynomial-time solution to the
conjugacy problem in $B_4$. Finally the {\em ultra summit set}
$$USS(x)=\{y\in SSS(x)\mid \mathbf c^M(y)=y
\mbox{ for some positive integer }M\}$$ was used by Gebhardt in
\cite{Ge} to propose a new algorithm together with experimental
data demonstrating the efficiency of his algorithm. Clearly
$$RSSS(x)\subset USS(x)\subset SSS(x)\subset SS(x),$$
and $RSSS(x)=USS(x)$ if $x$ is rigid.

On the other hand, fewer researches have been done to solve the
reducibility problem, perhaps due to lack of suitable tools.
Bestvina and Handel \cite{BH} invented the ``train track"
algorithm and solved the reduction problem for any surface
automorphism. Unfortunately, this algorithm is typically
exponential for the length of input described as automorphisms of
graphs. Bernardete, Nitecki and Guti\'{e}rrez \cite{BNG} showed
that a standard reduction system is preserved by cycling and
decycling and so for any reducible braid $x$, some braid in
$SSS(x)$ must have a standard reduction system and consequently
the reduction problem can be solve as soon as $SSS(x)$ is
generated. Humphries \cite{Hum} solved the problem of recognizing
split braids.

Recently there have two noticeable progresses. Ko and J.~Lee
\cite{KL}, Birman, Gebhardt and Gonz\'{a}lez-Meneses \cite{BGG1}
showed that some power of a pseudo-Anosov braid is rigid, up to
conjugacy and in fact any braid in the ultra summit set of the
power is rigid. E.~Lee and S.~Lee \cite{L05} showed that if the
outermost component of a reducible braid $x$ is simpler than the
whole braid $x$ up to conjugacy then any braid in $RSSS(x)$ has a
standard reduction system. The difficulty of using the $USS$ or
the $RSSS$ is that we do not know how long it takes to generate
one element in the set, not to mention the whole set.

Our contribution in this article is two folds. We now know how
fast we can obtain a rigid braid from a given pseudo-Anosov
$n$-braid by taking powers and iterated cyclings. In fact the
required power is at most $(n(n-1)/2)^3$ and the required number
of iterated cyclings is at most $2n!(n(n-1)/2)^3\ell(x)$. The
other contribution is a complete understanding of reducible braids
that are rigid. In fact, if a reducible braid is rigid and no
circles in its reduction system are standard, its conjugate by
some permutation braid must preserve a set of standard circles.

We will start by introducing these contribution as two main
theorems in the next section as well as a polynomial-time
algorithm for the reducibility problem. Proofs for the main
theorems will follow in the next couple of sections.

The authors wish to thank Joan Birman, Volker Gebhardt and
Juan Gonz\'{a}lez-Mensese for many helpful comments made toward the first
draft of this paper.

%All four
%invariant sets enjoy the property that if $a^{-1}ya\in P$ and
%$b^{-1}yb\in P$ for $y\in P$ and $a,b\in S_n$ then $(a\wedge
%b)^{-1}t(a\wedge b)\in P$ where $P$ denotes one of invariant sets.

%If $\inf(x)<\inf_c(x)$, then an
%iterated cycling on $x$ increases the infimum, that is,
%$\inf(x)<\inf(\mathbf c^N(x))$ for some positive integer $N$
%\cite{EM}. It is known that $N\le n(n-1)/2$ \cite{BKL}. Since
%$\sup(x)=-\inf(x^{-1})$ and $\tau(\mathbf{d}(x))=\mathbf
%{c}(x^{-1})^{-1}$, the statement obtained by replacing $\inf$,
%$\mathbf c$ by $\sup$, $\mathbf d$ and reversing inequalities also
%holds.

\section{Main theorems and an algorithm for the reducibility problem}
We introduce two main theorems whose proofs constitute the next
two sections. Then we present an algorithm for the reducibility
problem based on the main theorems. Throughout this paper we will
assume that $n\ge 3$ since every 2-braid is periodic.

\subsection{Main theorems}
The word length of the fundamental $n$-braid $\Delta$ is
$\frac{n(n-1)}{2}$ and will be denoted by $D$.

\begin{thm} \label{thm:1}
Let $x$ be a pseudo-Anosov $n$-braid. Then there are integers
$1\le L\le 2D$, $1\le M\le D^2$, and $1\le N\le n!\,\ell(x)LM$
such that the weighted form of $\mathbf c^N(y^L)$ is rigid for any
$y\in SSS(x^M)$.
\end{thm}

\begin{thm} \label{thm:2}
Let $x$ be a reducible, rigid $n$-braid. Then there exists a
permutation $n$-braid $t$ such that $t^{-1}xt$ is rigid and has at
least one orbit of standard reduction circles.
\end{thm}

\subsection{Reducibility Algorithm}\mbox{}\\
Input: An $n$-braid $x$ given as a word in the Artin generators\\
Output: The dynamical type of $x$, that is, whether $x$ is
periodic, pseudo-Anosov, or reducible.
\begin{enumerate}
    \item We first try to find a positive integer $M$ such that
    $\inf_c((x^M)^D)=D\inf_c{(x^M)}$ and $\sup_c((x^M)^D)=D\sup_c{(x^M)}$
    and choose $y \in SSS(x^M)$. By \cite{LL01}, such an $M$ must exist
    between 1 and $D^2$ and so we can obtain it via the loop:
    For $1\le j\le D^2$,
    we test whether $$\inf(\mathbf d^{D\ell(x)jD}\mathbf c^{D\ell(x)jD}(x^{jD}))
    =D\inf(\mathbf d^{D\ell(x)j}\mathbf c^{D\ell(x)j}(x^j))$$
    and $$\sup(\mathbf d^{D\ell(x)jD}\mathbf c^{D\ell(x)jD}(x^{jD}))
    =D\sup(\mathbf d^{D\ell(x)j}\mathbf c^{D\ell(x)j}(x^j))$$
    after computing necessary weighted forms, and then
    if the test answers affirmatively, return $M=j$;

    \item If $\mathbf d^{D\ell(x)M}\mathbf c^{D\ell(x)M}(x^M)=\Delta^{2k}$
    for some $k$, then conclude that $x$ is periodic and stop.
    Otherwise, set $y=(\mathbf d^{D\ell(x)M}\mathbf
    c^{D\ell(x)M}(x^M))^{2D}$;

    \item Test whether there exists an integer $1\le N\le n!\ell(y)$
    such that $\mathbf{c}^{N}(y)$ is rigid. If such an $N$ does not exist,
    then conclude that $x$ is reducible and stop. Or if $\mathbf{c}^{N}(y)$ is rigid,
    set $z=\mathbf{c}^{N}(y)$.

    \item Test whether there exists a permutation braid $t\in S_n$ such that
    $t^{-1}zt \in RSSS(z)$ and $t^{-1}zt$ has at least one orbit of standard
    reduction circles.
    If such a $t$ exists, conclude that $x$ is reducible.
    Otherwise, conclude that $x$ is pseudo-Anosov.
\end{enumerate}

We now explain why our algorithm works and analyze its complexity
in step by step.

In Step (1), notice that $\inf_c(\beta)=\inf(\mathbf
c^{\ell(\beta)D}(\beta))$ since if $\inf_c(\beta)>\inf(\beta)$,
$\inf(\mathbf c^{D}(\beta))>\inf(\beta)$ by \cite{BKL} and
$\inf_c(\beta)\le \inf(\beta)+\ell(\beta)$. Similarly,
$\sup_c(\beta)= \sup(\mathbf d^{\ell(\beta)D}(\beta))$. Thus
$\mathbf d^{\ell(\beta)D}\mathbf c^{\ell(\beta)D}(\beta)\in
SSS(\beta)$. The complexity of this step is dominated by
$\ell(x)jD^2$ times of cyclings and decyclings on $x^{jD}$ for
$j=1,\ldots,D^2$. Thus it can be estimated as $\mathcal
O(\ell(x)^3n^{21}\log n)$.

Step(2) is simple and its complexity can be dominated by other
steps.

In Step (3), the existence of such an $N$ for any pseudo-Anosov
braid and the upper bound for $N$ are a part of
Theorem~\ref{thm:1}. The complexity of this step is $\mathcal
O(n!\,\ell(y)\cdot \ell(y)^2n\log n)=\mathcal
O(n!\,(2D)^3\ell(x)^3n\log n)=\mathcal O(\ell(x)^3n!\,n^{19}\log n)$.

Step (4) is the most complicated.  One can check that $t^{-1}zt$
has at least one orbit of standard reduction circles by an easy
algorithm such as one given in \cite{L05} since a standard
reduction circle is preserved by each permutation braid that is a
factor in the weighted form of $t^{-1}zt$. Indeed for each pair
$(i,j)$ with $1\le i<j\le n$, let $L=\lfloor \frac
n{j-i+1}\rfloor$ and $(t^{-1}zt)^L=\Delta^u z_1z_2\cdots
z_{L\ell(z)}$ be the weighted form and $\hat z_k$ denote the
bijection on $\{1,\ldots,n\}$ corresponding to the permutation
braid $z_k$. We need to check whether $\{\hat z_1\hat z_2\cdots
\hat z_m(i), \hat z_1\hat z_2\cdots \hat z_m(i+1),\ldots, \hat
z_1\hat z_2\cdots\hat z_m(j)\}=\{i, i+1,\ldots,j\}$ for some $1\le
m\le L\ell(z)$ that is a multiple of $\ell(z)$ and the set $\{\hat
z_1\hat z_2\cdots \hat z_k(i), \hat z_1\hat z_2\cdots \hat
z_k(i+1),\ldots, \hat z_1\hat z_2\cdots\hat z_k(j)\}$ is consisted
of consecutive integers for each $1\le k\le m$. The complexity of
these two test is $\mathcal O(\ell(z)^2n\log n \cdot
\ell(z)n^6)=\mathcal O(\ell(x)^3n^{25}\log n)$.

A naive algorithm would be to perform these two tests for each
$t\in S_n$ and then the complexity of Step(4) is $\mathcal
O(\ell(x)^3n!\,n^{25}\log n)$. Since $t^{-1}zt \in RSSS(z)$ iff
$t^{-1}zt$ is rigid \cite{BGG1,KL}, conjugators preserving rigid
braids are closed under meet. Thus we may improve our algorithm in
practice.  Indeed find all minimal braids $t_1, t_2, \ldots, t_m$
for some $1\le m\le n-1$ such that $t_i^{-1}zt_i\in RSSS(z)$.
These $t_i$ can be found by a formula starting from generators.
For each $t_i^{-1}zt_i$, perform the test for the possession of a
standard circle. If nothing is found, inductively find the next
larger minimal conjugators preserving $RSSS(z)$ by starting from a
join of two minimal conjugators found in the previous steps. The
latter method would be much faster in the average case but the
complexity in the braid index $n$ is not clear yet. Thus the
over-all complexity is cubic in the canonical length of the input
braid.

\section{Proof of Theorem \ref{thm:1}}
\begin{lem}\label{lem:cwpower}
Suppose that $x$ is an $n$-braid such that $x\in SSS(x)$ and
$\inf(x^i)=i\inf(x)$, $\sup(x^i)=i\sup(x)$ for $i\ge 1$. If $x^K$ is
rigid for some $K\ge 1$ then $x$ itself is rigid.
\end{lem}
\begin{proof} Under the hypotheses, neither new $\Delta$'s can be
formed nor factors can be merged by taking powers. Thus
$\tau^{u(K-1)}(H(x))\prec H(x^K)$ and $T(x)\succ_R T(x^K)$.
Consequently $T(x^K)\lceil\tau^{uK}(H(x^K))$ implies
$T(x)\lceil\tau^{u}(H(x))$.
\end{proof}

\begin{pro}[\cite{BGG1,KL}]\label{thm:pacw}
Let $x$ be a pseudo-Anosov braid in $B_n$. Then $x^M$ is conjugate
to a rigid braid for some $1\le M\le (\frac{n(n-1)}{2})^2$.
\end{pro}

\begin{proof}
we refer to \cite{BGG1,KL} for the existence of such an $M$ and we
add a comment on the upper bound for $M$. By Theorem 4.3 in
\cite{LL01}, there exist a positive integer $M\leq
(\frac{n(n-1)}{2})^2$ and $y\in C(x)$ such that
$\inf((y^M)^i)=i\inf(y^M)$ and $\sup((y^M)^i)=i\sup(y^M)$ for
$i>0$. Let $z=y^M$. Since $z$ is pseudo-Anosov, $z^{M'}$ is
conjugate to a rigid braid for some $M'>0$. Hence, $z$ is
conjugate to a rigid braid by Lemma \ref{lem:cwpower}.
\end{proof}

A braid $x$ is {\em tame} if $\inf(x^i)=i\inf(x)$ and
$\sup(x^i)=i\sup(x)$ for all $i\ge 1$. For any $n$-braid $y$,
$y^M$ becomes tame for some $1\le M\le D^2$ by \cite{LL01}.

\begin{lem}\label{tame}
If $x, y\in SSS(x)$ for a tame braid $x$, then $y$ is also tame.
\end{lem}
\begin{proof}
It is clear that $\inf_c(y^i)\ge \inf(y^i)\ge i\inf(y)$ for all
$i\ge 1$. By the hypothesis, $x^i\in SSS(x^i)$ and so
$\inf_c(y^i)=\inf_c(x^i)=\inf(x^i)=i\inf(x)=i\inf(y)$ for all
$i\ge 1$. Thus $\inf(y^i)=i\inf(y)$. Similarly,
$\sup(y^i)=i\sup(y)$ for $i\ge 1$.
\end{proof}

\begin{lem}\label{tame to i-rigid}
Let $x$ be a tame braid and $i\ge 1$. Then $x^L$ is left (or
right) $i$-rigid for some $1\le L\le iD$. In particular, $x^{iD}$
is left and right $i$-rigid.
\end{lem}
\begin{proof}
For the simplicity of notations, assume that $\inf(x)=0$ and prove
the left version. Let $H_i(y)$ denote the product of the first $i$
factors in the weighted form of a tame braid $y$. Since $H_i(x^L)$
can not be strictly increasing without producing $\Delta$ for all
$L$ that increases from 1 to $iD$, $H_i(x^L)=H_i(x^{L+1})$ for
some $1\le L\le iD$. Then
$H_i(x^{M})=H_i(x^{M-L-1}H_i(x^{L+1}))=H_i(x^{M-L-1}H_i(x^{L}))=H_i(x^{M-1})$
for all $M \ge L+1$. Thus $H_i(x^{L})=H_i(x^{M})$ for all $M \ge
L$. In particular, $H_i(x^{M})=H_i(x^{2M})$ and so $x^M$ is
$i$-rigid for all $M \ge L$.
\end{proof}

\begin{lem}\label{i-rigid}
Suppose that $x$ is tame and $x\in SSS(x)$. Then $\mathbf
c^j(x^{iD})$ is left and right $i$-rigid for all $j\ge 0$.
\end{lem}
\begin{proof}
By Lemma~\ref{tame to i-rigid}, $x^{iD}$ is left and right
$i$-rigid. It is enough to show that $\mathbf c(x^{iD})=y^{iD}$
for some tame braid $y\in SSS(y)$. Then we ought to set
$y=\tau^{u}(H(x^{iD})^{-1})x\tau^{u}(H(x^{iD}))$ for
$u=\inf(x^{iD})$. Since $(x^{iD})^{-1}x(x^{iD})=x\in SSS(x)$,
$y\in SSS(x)=SSS(y)$. By Lemma~\ref{tame}, $y$ is also tame.
\end{proof}

\begin{lem}\label{2-rigid}
If $x$ is left 2-rigid,
$H(\tau^u(a^{-1})x)=\tau^u(H(\tau^u(a^{-1})x^2))$ for any $a\prec
H(x)$ where $u=\inf(x)$. The corresponding statement using right
versions also holds.
\end{lem}
\begin{proof}
Let $x=\Delta^u x_1\cdots x_k$ be the weighted form. Then $a\prec
x_1$. Since the 2-rigidity implies $$H(x_2\cdots x_k)=H(x_2\cdots x_k\tau^u(x_1\cdots
x_k)),$$ we have
$$H((a^{-1}x_1)x_2\cdots x_k)=H((a^{-1}x_1)x_2\cdots x_k\tau^u(x_1\cdots x_k)).$$
\end{proof}

\begin{thm}\label{cycling bound and WCW}
Let an $n$-braid $x$ be the $2D$-th power of a tame braid that is
in its super summit set. If $x$ is conjugate to a rigid braid,
then a rigid braid must be obtained from $x$ by at most
$n!\,\ell(x)$ iterated cyclings.
\end{thm}
\begin{proof}
Since we are assuming $n\ge 3$, $\ell(x)\ge 6$. It was proved in
Theorem 3.3 in \cite{KL} and Theorem 3.15 in \cite{BGG1} that if
$USS(x)$ contains at least one rigid braid, then every braid in
$USS(x)$ is rigid. Thus iterated cyclings on $x$ must produce a
rigid braid. Let $y=\mathbf c^N(x)$ be the rigid braid obtained
from $x$ by the minimal number of iterated cyclings. Since
$\inf(x)$ is even, we can assume $\inf(x)=0$ for the sake of
simplicity without affecting the conclusion.

Let $y=y_1y_2\cdots y_k$ be the weighted form. We will prove by
induction on $i\ge 1$ that for all $1\le i\le N$
$$\mathbf c^{N-i}(x)=a_{i}y_{[1-i]}z_{i}$$
for some permutation braid $a_{i}$ satisfying $y_{[2-k-i]}\cdots
y_{[-1-i]}y_{[-i]}\succ_R a_{i}\succneqq_R e$ and a positive braid
$z_{i}=y_{[2-k-i]}\cdots y_{[-1-i]}y_{[-i]}a_i^{-1}$ where $[m]$
denotes the integer between 1 and $k$ that equals $m$ mod $k$. Let
$t_i=H(\mathbf c^{N-i-1}(x))$ for $1\le i\le N$. Then $\mathbf
c^{N-i-1}(x)=t_i\mathbf c^{N-i}(x)t_i^{-1}$ and the properties
\begin{enumerate}
\item[(i)] $t_i\lceil \mathbf c^{N-i}(x)t_i^{-1}$; \item[(ii)]
$\mathbf c^{N-i}(x)\succ_R t_i$; \item[(iii)] $t_i\succneqq_R
T(\mathbf c^{N-i}(x))$
\end{enumerate}
are clear from the definition of cycling, except the fact that
$t_i\not= T(\mathbf c^i(x))$. If $t_i= T(\mathbf c^i(x))$, then
$\mathbf c^{i-1}(x)$ is already rigid and this violates the
minimality of $N$.

Let $a_1=t_1y_k^{-1}$. Then $a_1$ is a permutation braid such that
$y_1\cdots y_{k-1}\succ_R a_1 \succneqq_R e$ by (ii) and (iii).
Thus our claim is proved for $i=1$. Suppose that our claim holds
for $i$. Then $\mathbf
c^{N-i-1}(x)=t_ia_{i}y_{[1-i]}z_{i}t_i^{-1}$. We must have
$t_{i}a_{i}\succ_R y_{[-i]}$, otherwise the factor $y_{[-i]}$
splits into two parts in $\mathbf c^{N-i-1}(x)$ so that
$\ell(\mathbf c^{N-i-1}(x))>\ell(x)$ by (i) which is a contradiction.
Thus we can write $t_{i}a_{i}=a_{i+1}y_{[-i]}$ for some
permutation braid $a_{i+1}$ satisfying $a_{i+1}\prec t_{i}$ and
$\mathbf c^{N-i}(x)\succ_R a_{i+1}$. By the minimality of $N$,
$a_{i+1}\not= e$. Let $H_R(w)$ denote the right head of $w$. Then
\begin{eqnarray*}
H_R(\mathbf c^{N-i}(x))&=& H_R((\mathbf c^{N-i}(x))^2)\\
&=& H_R(a_{i}y_{[1-i]}z_{i}a_{i}y_{[1-i]}z_{i})\\
&=& H_R(a_{i}y_{[1-i]}\cdots y_{[-1-i]}y_{[-i]}y_{[1-i]}\cdots y_{[-1-i]}y_{[-i]}a_i^{-1})\\
&=& H_R(y_{[1-i]}\cdots y_{[-1-i]}y_{[-i]}a_i^{-1})\\
&=& H_R(y_{[1-i]}z_{i}).
\end{eqnarray*}
Here, the first equality holds since $\mathbf c^{N-i}(x)$ is right
1-rigid by Lemma~\ref{i-rigid} and the the fourth equality follows
from Lemma~\ref{2-rigid} since $\mathbf
c^{N+[-i]}(x)=y_{[1-i]}\cdots y_{[-1-i]}y_{[-i]}$ is right 2-rigid
by Lemma~\ref{i-rigid}. Thus
$$y_{[2-i]}\cdots y_{[-1-i]}y_{[-i]}\succ_R y_{[1-i]}z_{i}\succ_R a_{i+1}\succneqq_R e$$
and $\mathbf c^{N-i-i}(x)=a_{i+1}y_{[-i]}z_{i+1}$ for
$z_{i+1}=y_{[1-k-i]}\cdots y_{[-2-i]}y_{[-1-i]}a_{i+1}^{-1}$ and
this completes the induction.

By our claim, $\mathbf c^{N-i}(x)$ is completely determined by
choosing a nontrivial permutation braid $a_i$ satisfying
$a_i\prec_R H_R(y_{[2-k-i]}\cdots y_{[-1-i]}y_{[-i]})$. For each
$1\le i\le \ell(x)$, there are at most $n!$ such choices. Thus
$N\le n!\,\ell(x)$.
\end{proof}
We remark that the braid $x$ in Theorem~\ref{cycling bound and WCW} need not be
pseudo-Anosov.

\section{Proof of Theorem \ref{thm:2}}
Theorem~\ref{thm:2} is used in Step (4) of our algorithm in Section~2 and all braids
dealt in Step (4) has the canonical length $\ge 6$. In this section,
we will assume that the braid index is $\ge 3$ and the canonical length is $\ge 2$ unless stated otherwise in order to avoid any unnecessary nuisance.

Let $P$ denote one of conjugacy invariant sets $SS, SSS, USS,
RSSS$ and let $y\in P(x)$. If a nontrivial positive $n$-braid
$\gamma$ satisfies $\gamma^{-1}y\gamma \in P(x)$, $\gamma$ is
called a {\em $P$-conjugator} of $y$. A $P$-conjugator $\gamma$ of
$y$ is {\em minimal} if either $\gamma\prec \beta$ or
$\gamma\wedge \beta=e$ for each positive braid $\beta$ with
$\beta^{-1}y\beta \in P(x)$. In fact it is not hard to see that a
minimal $P$-conjugator satisfies $\gamma\prec
\tau^{\inf(y)}(H(y))$ or $\gamma\prec T(y)^*$ or both (For
example, see \cite{BGG2,KL}). A conjugator $\gamma$ satisfying
$\gamma\prec \tau^{\inf(y)}(H(y))$ (or $\gamma\prec T(y)^*$,
respectively) will be called a {\em cut-head} (or {\em add-tail})
conjugator. In particular, if $y$ is rigid and $\gamma$ is its
$P$-conjugator then it can not be both cut-head and add-tail since
$T(y)\lceil \tau^{\inf(y)}(H(y))$, that is, $T(y)^*\wedge
\tau^{\inf(y)}(H(y))=e$. If $\gamma$ is a $USS$-conjugator of a
rigid braid $y$, it is also a $RSSS$-conjugator and
$\gamma^{-1}y\gamma$ is rigid (For example, see \cite{BGG2,KL}).
We note that if $\gamma$ is an add-tail conjugator of $y$, then
$\gamma$ is a cut-head conjugator of $y^{-1}$ (For example, see
\cite{BGG2,KL}).

Suppose that $\gamma$ is a product of cut-head conjugators of $y$,
that is, $\gamma=\gamma_1\cdots \gamma_j$ such that $\gamma_j$ is a
cut head conjugator of
$$(\gamma_1\cdots
\gamma_{j-1})^{-1}y(\gamma_1\cdots \gamma_{j-1})$$ where
$\gamma_0=e$. By the definition of cut head conjugators,
$\gamma\prec \tau^{\inf(y^i)}(\Delta^{-\inf(y^i)}y^i)$ for some
$1\le i\le j$. Similarly, if $\gamma$ is a product of add-tail
conjugators of $y$, then $\gamma\prec
\tau^{\inf((y^{-1})^{i})}(\Delta^{-\inf((y^{-1})^{i})}(y^{-1})^{i})$
for some $1\le i\le j$.
\begin{lem}
\label{subword} Suppose that an $n$-braid $x$ is rigid
and $\ell(x)\ge 2$. If $\gamma$ is a positive $n$-braid such that
$\gamma^{-1}x\gamma$ is rigid and $\ell(\gamma)\ge 2$
then there is a positive braid $\beta$ such that $\beta^{-1}x\beta$
is also rigid, $\ell(\beta)\ge 2$, and moreover
$\beta$ is a left subword of either
$\gamma\wedge\tau^{\inf(x^i)}(\Delta^{-\inf(x^i)}x^i)$ or
$\gamma\wedge\tau^{\inf((x^{-1})^i)}(\Delta^{-\inf((x^{-1})^i)}x^i)$
for some $i\ge 1$.
\end{lem}

\begin{proof}
If $h$ is a cut-head $RSSS$-conjugator of $x$ and $t$ is an
add-tail $RSSS$-conjugator of $h^{-1}xh$, then it is shown in
Proposition 3.23 in \cite{BGG2} that $ht$ is a permutation braid
and there are an add-tail $RSSS$-conjugator $t'$ of $x$ and a
cut-head $RSSS$-conjugator $h'$ of $t'^{-1}xt'$ such that
$ht=t'h'$. Thus we can write $\gamma=HT=T'H'$ such that $H$ (or
$H'$, respectively) is a product of cut-head $RSSS$-conjugators of
$x$ (or $T'^{-1}xT'$) and $T$ (or $T'$, respectively) is a product
of add-tail $RSSS$-conjugators of $H^{-1}xH$ (or $x$). If
$\ell(\gamma)\ge 2$, at least one of $\ell(H)$ or $\ell(T')$ is
$\ge 2$. Then the conclusion follows from the remark right before
this lemma.
\end{proof}

\begin{pro}\label{cycling and minimal conjuagtor}
Suppose that an $n$-braid $x$ is rigid and $\alpha^{-1}x\alpha\in
SSS(x)$ for a permutation $n$-braid $\alpha$. Let $\mu$ be minimal
among all $\beta$'s satisfying that $\alpha \prec \beta$ and
$\beta^{-1}x\beta \in RSSS(x)$. Then $\mathbf
c^N(\alpha^{-1}x\alpha)=\mu^{-1}x\mu$ for some $N\ge 0$.
\end{pro}
\begin{proof}
Since $x$ is rigid, $\mathbf c^{2i\ell(x)}(x)=x$ for $i\ge 1$. Note
that there exists $\alpha'$ such that $\mathbf
c^{2\ell(x)}(\alpha^{-1}x\alpha)=\alpha'^{-1}x\alpha'$ and
$\alpha\prec \alpha'$. Let $\mathbf
c^{2i\ell(x)}(\alpha^{-1}x\alpha)=\alpha_i^{-1}x\alpha_i$ where
$\alpha_i\prec \alpha_{i+1}$ and $\alpha_0=\alpha$. Then
$\alpha_j^{-1}x\alpha_j\in RSSS(x)$ for some $j>0$ since $x$ is
rigid. Let $\mu$ be minimal among all $\beta$'s satisfying that
$\alpha \prec \beta$ and $\beta^{-1}x\beta \in RSSS(x)$. Since
$\alpha\prec \mu$, $\alpha_i\prec \mu$ for $i\ge 0$ by Corollary 2.2
in \cite{Ge}. Thus $\alpha_j\prec \mu$ and so $\alpha_j= \mu$ by
the minimality of $\mu$.
\end{proof}

\begin{figure}[h]
\center
\includegraphics{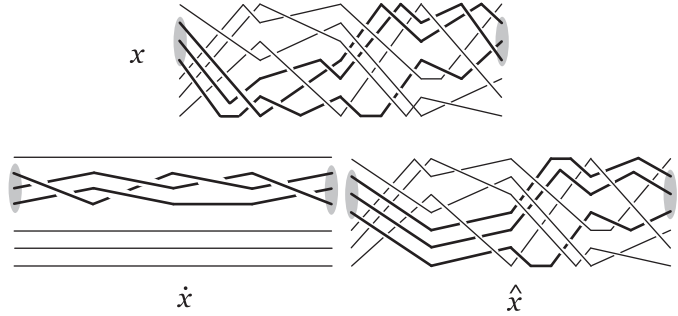}
\caption{The decomposition of a reducible
braid}\label{fig:decomposition}
\end{figure}

If a reducible $n$-braid $x$ has an orbit of standard reduction
circles, we may assume that $x$ preserves a standard circle $C$ by
replacing $x$ by $x^j$ for some $1\le j\le n$ if necessary. $C$
spans a tube that does not intersect any strand of $x$. Then the
standard reduction circle $C$ uniquely determines a decomposition
$x=\dot x\hat x=\hat x\dot x$ where all strands lying outside $C$
and the tube spanned by $C$ form a trivial braid in $\dot x$ and on
the other hand all strands lying inside $C$ form a trivial braid in
$\hat x$ as in Figure~\ref{fig:decomposition}.

\begin{figure}[h]
\center
\includegraphics{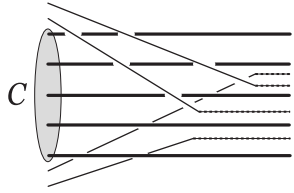}
\caption{A typical destroyer of $C$}\label{fig:destroyer}
\end{figure}

For a standard circle $C$, a permutation $n$-braid $\beta$ is called a {\em destroyer} of $C$ if
$\alpha(C)$ is not standard for all $e\precneqq \alpha \prec \beta$. See
Figure~\ref{fig:destroyer}.
Suppose the standard circle $C$ contains $\ell$ consecutive punctures from the $k$-th
puncture where we must have either $k\ge 2$ or $k+\ell\le n-1$ to
make $C$ a reduction circle. Let $\beta$ be a destroyer of $C$, then $\beta(C)$ is not standard,
$\beta$ should ends with a nontrivial permutation $n$-braid
$\delta$, i.e. $\beta\succ_R \delta$ such that either
$\delta\precneqq\sigma_{k-1+s}\cdots \sigma_{k+\ell+s}$ or
$\delta\precneqq\sigma_{k+\ell+s}\cdots \sigma_{k-1+s}$ for some
integer $s$ satisfying $k-1+s\ge 1$ and $k+\ell+s\le n-1$. Furthermore,
since every left subword of $\beta$ must destroy the standard circle
$C$, we have that
$\dot\beta=e$, there are no crossings among outer strands, and no outer strands
pass through all of inner strands at once.
Thus a typical destroyer of $C$ should look like one given in
Figure~\ref{fig:destroyer} where $\dot\beta$ is drawn by thicker strands.

\begin{lem}\label{destroyer} Given a standard circle $C$ and
permutation braids $t$, $y$, and $s$, suppose that $ty$ sends $C$
to a standard circle, $(ty)\lceil s$, and $ty, ys\prec \Delta$. If
$t$ is a destroyer of $C$ then $s$ is also a destroyer of $ty(C)$.
\end{lem}
\begin{proof}
Let $C'=ty(C)$. Assume that $s$ is not a destroyer of $C'$, that is,
$s'(C')$ is standard for some $e\precneqq s' \prec s$. Suppose $s'$ is minimal
among such braids. We have either
$s'(C')=C'$ or $s'(C')$ is another standard circle.
If $s'(C')=C'$ then $s'$ has one of
two types in Figure~\ref{fig:type of s}.
\begin{figure}[h]
\center Type I \includegraphics{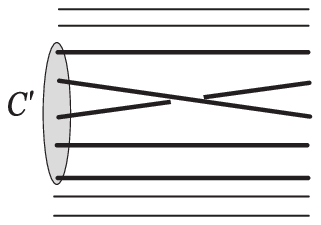}
\qquad\includegraphics{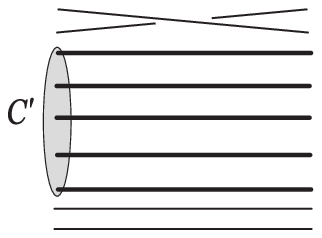} Type II
\caption{Possible $s'$ when $s'(C')=C'$}\label{fig:type of s}
\end{figure}
If $s'$ is of Type I, then $t$ must has crossings of strings in the tube
spanned by $C$ to satisfy $(ty)\lceil s$, but this contradicts the assumption that $t$ is a destroyer of $C$. If $s'$ is of Type II, then there are three possibilities for $t$ to satisfy $(ty)\lceil s$
as in Figure~\ref{fig:cases}.
\begin{figure}[h]
\center
\includegraphics{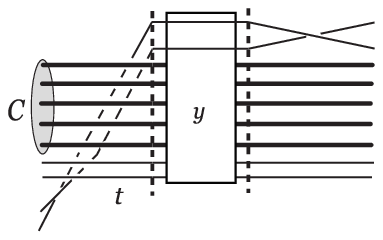}
\quad
\includegraphics{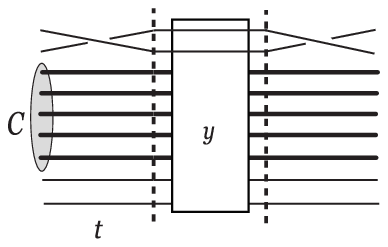}
\quad
\includegraphics{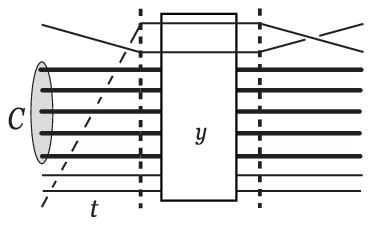}
\caption{Possible $t$ when $s'$ is of type II}
\label{fig:cases}
\end{figure}
In all cases, $t$ must have a subword preserving a circle $C$ but
this contradicts the assumption that $t$ is a destroyer of $C$.

If $s'(C')$ is another standard circle then $s'$ must be of the
type in Figure~\ref{fig:case IV}. Thus $t$ should look like
Figure~\ref{fig:case IV} to satisfy $(ty)\lceil s$. Again $t$ must
have a subword preserving a circle $C$ but this contradicts by the
assumption that $t$ is a destroyer of $C$. Consequently, $s$ is a
destroyer of $C'$.

\begin{figure}[h]
\center
\includegraphics[scale=1.2]{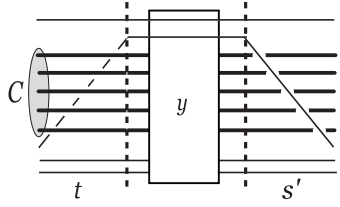}
\caption{Possible $s'$ and $t$ when $s'(C')$ is another standard circle}
\label{fig:case IV}
\end{figure}
\end{proof}

\begin{lem}
\label{forbidden head}  Suppose that a reducible $n$-braid $x$ is
rigid and has a standard circle $C$. If $\beta$ is a permutation
$n$-braid such that $\beta^{-1}x\beta$ is rigid and $\alpha(C)$ is
not standard for any $e\precneqq \alpha\prec \beta$ such that
$\alpha^{-1}x\alpha$ is rigid, then $\beta$ is a destroyer of $C$.
\end{lem}
\begin{proof}
For the simplicity, assume $\inf(x)=0$.
Suppose that $\beta$ is not a destroyer of $C$. Then $\alpha(C)$ is
standard for some $e\precneqq \alpha\prec \beta$. We choose a maximal
such $\alpha$ so that $t=\alpha^{-1}\beta$ is a destroyer of $\alpha(C)$.
We will show that $\alpha^{-1}x\alpha\in SSS(x)$. Then
$\mathbf c^N(\alpha^{-1}x\alpha)=\mu^{-1}x\mu$ for some $\alpha\prec\mu$
and $\mu$ is minimal among all braids $\gamma$ such that
$\alpha\prec \gamma$ and $\gamma^{-1}x\gamma$ is rigid by Proposition
\ref{cycling and minimal conjuagtor}. Thus $e\precneqq \mu\prec
\beta$. Moreover since $\mu^{-1}x\mu$ is an iterated cycling of $\alpha^{-1}x\alpha$
that fixes a standard circle $\alpha(C)$, it fixes a standard circle $\mu(C)$ by \cite{BNG}.
However this contradicts the hypotheses.

Suppose $\alpha^{-1}x\alpha\notin SSS(x)$. Then $\ell(\alpha^{-1}x\alpha)=\ell(x)+1$.
Let $k=\ell(x)$ and $\alpha^{-1}x\alpha=z_1\cdots z_kz_{k+1}$ be the weighted form.
Then the weighted form of $\beta^{-1}x\beta=t^{-1}z_1\cdots z_kz_{k+1}t$
can be written
$$(t_1^{-1}z_1 t_2)\cdots
(t_{k-1}^{-1}z_{k-1}t_{k})(t_{k}^{-1}z_{k}z_{k+1}t_1)$$
where $t_1=t$ and $t_i\prec z_i$ for $1\le i\le k$.
Since $t_1=t$ is a destroyer of $\alpha(C)$ by the choice of $\alpha$,
$t_i$ is a destroyer of $C_i$ for $1\le i\le k$, where $C_1=\alpha(C)$ and
$C_i=z_1\cdots z_{i-1}(\alpha(C))$ for $2\le i\le k+1$ by Lemma
\ref{destroyer}. Since $z_{k}\lceil z_{k+1}$ and
$t_{k}$ is a destroyer of $C_{k}$, $z_{k+1}$ is a
destroyer of $C_{k+1}$ but this contradicts the fact that
$z_{k+1}(C_{k+1})$ is standard. Hence $\alpha^{-1}x\alpha\in SSS(x)$.
\end{proof}

\begin{thm}
\label{reducible}  Suppose that a reducible $n$-braid $x$ is
rigid and has an orbit of standard circles starting
with a standard circle $C$. If $\gamma$ is a positive $n$-braid such
that $\gamma^{-1}x\gamma$ is rigid and $\beta(C)$ is
not standard for any $e\precneqq \beta\prec \gamma$ such that
$\beta^{-1}x\beta$ is rigid, then $\ell(\gamma)\le 1$
\end{thm}
\begin{proof}

For the sake of simplicity, we assume that $\inf(x)=0$. If the
conclusion holds for a power of $x$, so does it for $x$ itself.
Thus we can further assume that $x(C)=C$ by replacing $x$ by its
power if necessary. Recall the notations for inner braids and
outer braids with respect to the standard curve $C$.

Suppose $\ell(\gamma)\ge 2$. By Lemma~\ref{subword}, we may assume
that $\gamma\prec x^i$ or $\gamma\prec (x^{-1})^i$ for some $i\ge
1$. Since $x^i(C)=C=(x^{-1})^i(C)$, we may work on either case. So
we assume $\gamma\prec x^i$. Let $\gamma_1\gamma_2$ be the first two
factors in the weighted form of $\gamma$.
By Lemma~\ref{forbidden head}, $\gamma_1$ must be a destroyer of $C$
as in Figure~\ref{fig:destroyer}.

Since $\gamma_1\gamma_2$ is left weighted, $\gamma_2$ must start
with one of two types of crossings given in Figure \ref{gamma12}. A
crossing of type I is formed by an inner strand and an outer stand.
On the other hand, a crossing of type II is formed by two outer strands
and is located between two inner strands. But we will show that both cases
are impossible.

\begin{figure}[h]\hbox{%
\parbox{.45\textwidth}{\center
\includegraphics{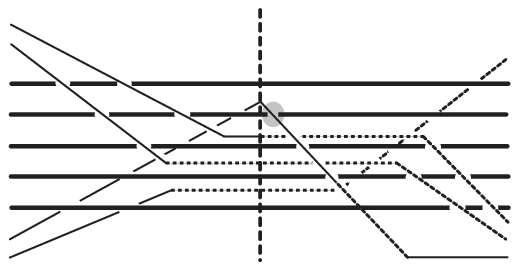}\\
Type I}
\hspace{.1\textwidth}
\parbox{.45\textwidth}{\center
\includegraphics{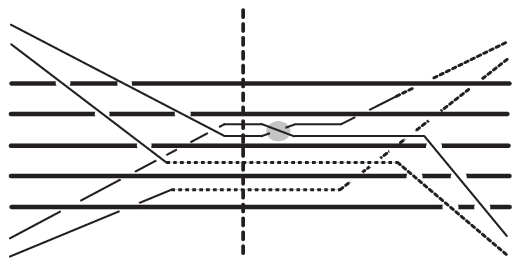}\\
Type II}}
\caption{Possible $\gamma_1\gamma_2$}
\label{gamma12}
\end{figure}

For any reducible braid $y$ fixing a standard circle $C$, consider the two-component link
$K_1\cup K_2$ obtained from $y$ by either one of the following ways:
\begin{enumerate}
\item[(i)] $K_1$ is obtained by the plat closing of any two inner strands and $K_2$ is obtained by Markov's closing of any outer strand.
\item[(ii)] $K_1$ is obtained by the plat construction of any two inner strands and $K_2$ is obtained by the plat closing of any two outer strands.
\end{enumerate}
It is clear that the link $K_1\cup K_2$ always splits, that is, there is an embedded 2-sphere separating two components.
However two components of the link obtained from $\gamma_1\gamma_2$ of type I via the construction (i) has the linking number 1.
Also two components of the link obtained from $\gamma_1\gamma_2$ of type II via the construction (ii) has the linking number 2. See Figure~\ref{gamma_close}.
\begin{figure}[h]\hbox{%
\parbox{.45\textwidth}{\center
\includegraphics{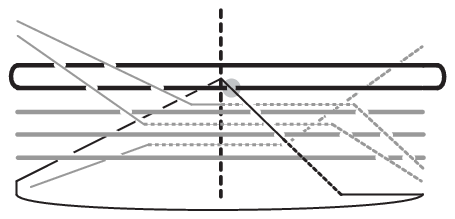}\\
Type I}
\hspace{.1\textwidth}
\parbox{.45\textwidth}{\center
\includegraphics{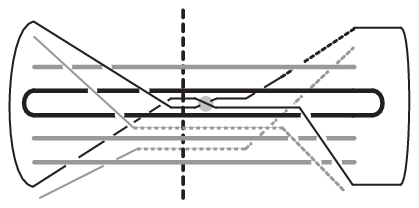}\\
Type II}}
\caption{Two-component links obtained from $\gamma_1\gamma_2$}
\label{gamma_close}
\end{figure}

Since $\gamma_1\gamma_2\prec\gamma\prec x^i$ and $x^i$ is a positive braid, this nontrivial linking can not be undone in $x^i$. This completes the proof.
\end{proof}

\begin{cor} Let $x$ be a reducible $n$-braid.
Suppose that $x$ is rigid and $x$ has no standard
reduction circle. If $\gamma^{-1}x\gamma$ is rigid
and has a standard reduction circle $C$ for $\gamma\in B_n^+$ with
$\ell(\gamma)>1$, then there exists $e\precneqq\beta\prec \gamma$
such that $\beta^{-1}x\beta$ is rigid and has at
least one standard reduction circle.
\end{cor}
\begin{proof}
Let $y=\gamma^{-1}x\gamma$. Then $x=\gamma y\gamma^{-1}$. Since
$\Delta^i(C)$ is a standard circle for all $i$ and
$(\gamma\wedge_R \dot \Delta)^{-1}(C)$ is a standard circle, we
can assume that $\inf(\gamma)=0$ and $\gamma\wedge_R \dot
\Delta=e$. Since $x=\gamma y\gamma^{-1}$,
$x=\Delta^{\ell(\gamma)}{\gamma^*}^{-1}
y\gamma^*\Delta^{-\ell(\gamma)}$, where
$\gamma\gamma^*=\Delta^{\ell(\gamma)}$. Let
$\gamma'=\tau^{\ell(\gamma)}(\gamma^*)$ and
$y'=\tau^{\ell(\gamma)}(y)$. If $\ell(\gamma')>1$ then there
exists $e\precneqq\beta'\prec \gamma'$ such that
$\beta'^{-1}y'\beta'$ is rigid and $\beta'(C)$ is a
standard circle by Theorem \ref{reducible}. Thus there exists
$\gamma \succ_R \beta'' \succneqq_R e$ such that
$\beta''y\beta''^{-1}$ has a standard circle $C$. Thus
$\beta^{-1}x\beta$ is rigid and has at least one
 standard reduction circle, where
$\beta=\gamma\beta''^{-1}$.
\end{proof}

By inductively applying Step 4 of the algorithm in Section~2,
one can find a whole reduction system for a reducible braid that is rigid after taking a power and iterated cyclings. The class of
these reducible braids is disjoint from the class of reducible
braids considered in \cite{L05}. In terms of an (outmost) orbit of standard
circles fixed by a reducible braid $x$, reducible braids considered
in \cite{L05} roughly satisfy $\sup(\dot x)>\sup(\hat x)$ and this
property is not inherited to $\dot x$ or $\hat x$, while rigid reducible braids satisfy $\sup(\dot x)=\sup(\hat x)$ and
both $\dot x$ and $\hat x$ are also rigid. The class
of reducible braids satisfying $\sup(\dot x)<\sup(\hat x)$ will also
need an attention.

\end{document}